\documentclass[a4paper,12pt]{article}
\usepackage{xspace}

\newcommand\NN{\mbox{I\kern-.21em N}}
\newcommand\RE{{\mbox{\rm I\kern-.21em R}}}   
\newcommand\ZZ{{\mbox{\sf Z\kern-.45em Z}}} 
\def\vv{\kern.344em{\rule[.18ex]{.075em}{1.32ex}}\kern-.344em}
\def\CX{{\mbox{\rm \vv C}}}
\newcommand\be{\beta}
\newcommand\gm{\gamma}
\newcommand\Gm{\Gamma}
\newcommand\dl{\delta}
\newcommand\eps{\varepsilon} \newcommand\f{\varphi}

\newcommand\Om{\Omega}
\newcommand\om{\omega}
\newcommand\lm{\lambda}
\newcommand\lmn{\lambda_n}
\newcommand\Lm{\Lambda}
\newcommand \tl {\tilde}
\newcommand\ds{\displaystyle}
\newcommand\iy{\infty}
\newcommand \underset [2]  { \mathop{#2}\limits_{#1}^{} }

\newcommand\ie{{inequality}}
\newcommand\lT{\ensuremath{L^2(0,T)}\xspace}
\newcommand\liy{\ensuremath{L^2(0,\iy)}\xspace}
\newcommand\op{operator\xspace}
\newcommand\uds{uniformly discrete set}
\newcommand\rds{relatively uniformly discrete sequence\xspace}

\newcommand\E{\ensuremath{{\mathcal E}}\xspace}

\newcommand\Lb{{${\mathcal L}$-basis}\xspace}
\newcommand\Rb{Riesz basis\xspace}
\newcommand\proof{\noindent {\sc Proof:}\qquad}
\newcommand\on[1]{\mathrm{#1}\,}
\newcommand\Np{\mathcal N^{(p)} }
\newcommand\p{^{(p)}}
\newfont{\Blackboard}{msbm10 scaled 1200}

\newfont{\roma}{cmr10 scaled 1200}
\newfont{\sBlackboard}{msbm10 scaled 900}
\newcommand{\sbl}[1]{\mbox{\sBlackboard #1}}

\newcommand{\szline} {{\sbl Z}}

\newcommand{\srline}  {{\sbl R}}

\newcommand{\KKK}{{\bf K}}

\renewcommand\d{\partial}
\renewcommand\a{\alpha}

\renewcommand\L{\mathcal L}
\renewcommand\iff{if and only if }

\renewcommand\d{\partial}

\newtheorem{definition}{\bf Definition}
\newtheorem{theorem}{\bf Theorem}
\newtheorem{remark}{\bf Remark}
\newtheorem{proposition}{\bf Proposition}
\newtheorem{lemma}{\bf Lemma}

\newcommand\qed{\hfill$\quad${\rule{3mm}{3mm}}\medskip\\}
\newcommand\LB[1]{\label{#1}}
\newcommand\BE[2]{\begin{#1} #2 \end{#1}}
\newcommand\EQ[2]{\BE{equation}{\LB{#1} #2}}
\newcommand\ARR[2]{\BE{array}{{#1} #2}}
\newcommand\EQA[3]{\EQ{#1}{\ARR{#2}{#3}}}

\def\l{\left} \def\r{\right}

\begin{document}
\begin{center}
{\LARGE \bf
Exponential Riesz bases of subspaces and divided differences
}
 \footnote {
 This work was partially supported by the Australian Research Council
(grant \# A00000723)
and by the Russian Basic Research Foundation
 (grant \# 99-01-00744).
  }

\vskip.3cm

{\large S.A. Avdonin
   \footnote{Department of Applied Mathematics and Control,
   St.Petersburg State University,
   Bibliotechnaya sq.~2,
   198904 St.Petersburg, Russia and
Department Mathematics and Statistics,
The Flinders University of South Australia,
 GPO Box 2100, Adelaide  SA 5001, Australia;
email: avdonin@ist.flinders.edu.au},
\
S.A.Ivanov
\footnote{Russian Center of Laser Physics,
St.Petersburg State University, Ul'yanovskaya 1,
198904 St.Petersburg, Russia; e-mail: sergei.ivanov@pobox.spbu.ru}
}
\end{center}

\begin{quotation}
\centerline {Abstract}
\noindent
Linear combinations of exponentials $e^{i{\lambda}_kt}$ in the
case where the distance between some points ${\lambda}_k$ tends to zero
are studied.
 D. Ullrich \cite{Ullrich} has proved the basis property of the divided
differences of exponentials in the case when
$\left\{\lambda_k\right\}=\bigcup \Lambda^{(n)}$ and
 the groups $\Lambda^{(n)}$
consist of equal number of points all of them are close enough to $n,\, n\in
{\bf Z}.$
 We have generalized this result
 for groups with arbitrary number of close points
 and obtained
a full description of Riesz bases of exponential divided differences.

\end{quotation}

\section{Introduction}\LB{n0}
Families of `nonharmonic' exponentials $\left\{e^{i\lambda_kt}\right\}$
appear in various fields of mathematics such as the theory of
nonselfadjoint operators (Sz.-Nagy--Foias model), the Regge problem for
resonance scattering, the theory of linear initial boundary value problems
for partial differential equations, control theory for distributed
parameter systems, and signal processing.  One of the central problems
arising in all of these applications is the question of the Riesz basis
property of an exponential family. In the space $L^2(0,T)$ this problem was
 considered for the first time in the classical work of R.~Paley and
N.~Wiener  \cite{PW}, and since then has motivated a great deal of work by
many mathematicians; a number of references are given in \cite{KNP},
\cite{Young}  and \cite{AI95}.  The problem was ultimately given a complete
solution \cite {Pavlov79}, \cite {KNP}, \cite {Minkin} on the basis of an approach
suggested by B.~Pavlov.

The main result in this direction can be formulated as follows
\cite{Pavlov79}.
\begin{theorem}\LB{Pavlov}
Let $\Lm:= \{\lambda_k | k\in \ZZ\}$ be a countable set of the complex
plane. The family $\left\{e^{i\lambda_kt}\right\}$ forms a Riesz basis in
$L^2(0,T)$ if and only if the following conditions are satisfied:

 (i) $\Lm$ lies in a strip parallel to the real axis,
$$
\sup_{k\in \ZZ}|\Im \lm_k| < \infty\,,
$$

and  is  uniformly discrete (or
separated), i.e.
\EQ{discr}{
\delta(\Lm):=\inf_{k\neq n}|\lm_k - \lm_n |>0\,;
}

(ii) there exists an entire function $F$ of
exponential type with indicator diagram of width $T$
and zero set $\Lm$
{\rm (the generating function of the family
$\left\{e^{i\lambda_kt}\right\}$ on the interval $(0,T)$)}
such that, for some real $h$, the function
$\left|F(x+ih)\right|^2$ satisfies
the Helson--Szeg\"{o}\ condition: functions
$
u, v \in L^{\infty}(\RE),
{\Vert v \Vert}_{L^{\infty}(\srline)} < \pi/2
$
 may be found such that
\EQ{HS}{
|F(x+ih)|^2=\exp\{u(x)+\tl v(x)\}
}
\end{theorem}

Here the map $v\mapsto \tilde{v}$
denotes the Hilbert transform
for bounded functions:
$$
 \tl v(x)=\frac 1\pi
 \int_{-\infty}^\infty v(t)
     \left \{ \frac 1{x-t}-\frac{t}{t^2+1} \right \} dt.
$$

\begin{remark} Note that (i) is equivalent to Riesz basis property of \E
in its span in $L^2(0,\iy)$ and (ii) is a criterion
that the orthoprojector $P_T$ from this span into $L^2(0,T)$ is an
isomorphism (bounded and boundedly invertible operator).
\end{remark}

It is  well known that the Helson--Szeg\"{o}\ condition
 is equivalent to
 the Muckenhoupt condition $(A_2)$:
$$
\sup_{I\in {\cal J}}\Bigl\{
\frac{1}{\left|I\right|}
\int_I\left|F(x+ih)\right|^2\,dx\,\,\frac{1}{\left|I\right|}
\int_I\left|F(x+ih)\right|^{-2}\,dx\Bigr\}\,\, < \,\,\infty,
$$
where ${\cal J}$ is the set of all intervals of the real axis.

The  notion of the generating function mentioned above
plays a central role in the modern theory of
nonharmonic
Fourier series \cite{KNP,AI95}.
This notion plays also an important role in the
theory of exponential bases in Sobolev spaces
(see \cite{AI00, IK, LyubSeip00}).
It is possible to write the explicit expression for this function
$$
F(z)=\lim_{R\to\infty}\prod_{|\lambda_k|\le R}(1-\frac{z_k}{\lambda_k})
$$
(we  replace  the  term  $(1-\lambda_k^{-1}z)$  by  $z$ if
$\lambda_k=0$).

The theory of nonharmonic Fourier series was successfully applied to control
problems for distributed parameter systems and formed the base of the
powerful method of moments (\cite{Butkovsky,R78,AI95}).
Recent investigations
into new classes of distributed systems such as hybrid systems,
structurally damped systems  have raised a number of new difficult
problems in the theory of exponential families (see, e.g.
\cite{HZ,MZ,JTZ}).
One of them is connected with
the properties of the family $\E=\left\{e^{i\lambda_kt}\right\}$
in the case when the set $\Lm$ does not satisfy the separation condition
(\ref{discr}), and therefore $\E$ does not form a Riesz basis in its span
in $\lT$ for any $T>0$.

Properties of such families in \liy have been studied for the first time in
the paper of V. Vasyunin \cite{Vasyunin78} (see also \cite[Lec. IX]{Nik}).
In the case when $\Lm$ is a finite
union of separated sets, a natural way to represent $\Lm$ as a set of
groups $\Lm^{(p)}$ of close points was suggested.  The subspaces spanned on
the corresponding exponentials form a \Rb. This means that there exists an
isomorphism mapping these subspaces into orthogonal ones. This fact
together with Pavlov's result on
 the orthoprojector $P_T$ (see Remark above)
gives a criterion of the \Rb property of  subspaces
of exponentials in \lT: the generating function have to satisfy  the
Helson--Szeg\"{o}\ condition (\ref {HS}).
Note that for the particular case when the
generating function is a sine type function (see definition in
\cite{Levin61}, \cite{KNP}, \cite{AI95}),  theorem of such a kind was
proved by Levin  \cite{Levin61}.

Thus,  Vasyunin's result and  Pavlov's geometrical approach give us
description of exponential Riesz bases of subspaces. If we do have a \Rb
of subspaces, clearly, we can choose an orthonormalized basis in each
subspace and obtain a \Rb of {\sl elements}. However, this way is not
convenient in applications when we need more explicit formulae. It is
important to obtain description of Riesz bases of elements which are `simple
and natural' linear combinations of exponentials.

The first result in this direction was obtained by D.~Ullrich \cite{Ullrich}
who considered sets $\Lm$ of the form $\Lm=\bigcup_{n\in \szline}\Lm^{(n)}$,
where subsets $\Lm^{(n)}$ consist of equal number (say, N)  real points
$\lm_1^{(n)},\ldots,\lm_N^{(n)}$ close to $n$, i.e., $| \lm_j^{(n)}-n|
< \eps$
for all $j$ and $n$.
 He proved that for sufficiently small $\eps >0$ (no estimate of $\eps$
was given)
the family of particular linear combinations of exponentials
$e^{i\lambda_kt}$ ---  the so--called {\em divided differences} constructed
by subsets $\Lm^{(n)}$ (see Definition \ref{gdd_def} in subsection
 \ref{n13}) ---
forms a Riesz basis in $L^2(0,2\pi N)$. Such functions arise in numerical
analysis \cite{Shilov}, and the divided difference of
$e^{i\mu t}$, $e^{i\lm t}$ of the first order is
$(e^{i\mu t}-e^{i\lm t})/(\mu-\lm)$. In a sense, the Ullrich result may be
considered as a perturbation theorem for the basis family
$\l\{e^{int}, te^{int}, \dots, t^{N-1}e^{int}\r\}, n\in {\bf Z}$.

The conditions of this theorem are rather restrictive and
 it can not be applied to
some problems arising in control theory (see, e.g. \cite{BKL,CZ,JTZ,LZ,MZ}).

 In the present paper we generalize Ullrich's result in several directions:
the set $\Lm$ is allowed to be  complex, subsets $\Lm^{(n)}$ are allowed to
 contain an arbitrary
number of points, which are not necessarily `very' close to each other (and,
moreover, to some integer).

Actually, we give a full description of
Riesz bases of
exponential divided differences and generalized divided differences
(the last ones appear in the case of
multiple points $\lm_n$). To be more specific, we take a {\it sequence}
$\Lm$ which is `a union' of a finite number of separated sets.
 Following Vasyunin we decompose $\Lm$ into groups
$\Lm^{(p)}$, then choose for each group the family of the
generalized divided differences (GDD) and prove that these functions
form a \Rb in \lT if the generating function of the
exponential family satisfies the Helson--Szeg\"{o}\ condition (\ref
{HS}). To prove that we show that GDD for points
$\lm_1, \dots, \lm_N$ lying in a fixed ball form `a uniform basis', i.e.
the basis constants do not depend on the positions of $\lm_j$ in the ball.
Along with that,  GDD depend on parameters analytically. Thus, this family
is a natural basis for the situation when exponentials $e^{i\lm t},
\lm \in \Lm$ do not
form even uniformly minimal family. For the particular case $\Lm=\{n\alpha\}
\cup \{n\beta\}, n \in {\bf Z},$ appearing in a problem of simultaneous
control,
this scheme was realized in \cite{AT}.

In the case when $\Lm$ is not a finite union of separated sets, we
present a negative result: for some ordering of $\Lm$,
GDDs do not form  a uniformly minimal family.
\begin{remark}
After this paper has been written, a result on Riesz bases of exponential
DD  in their span in $L^2(0,T)$ for
large enough  $T$ has been announced in
\cite{MR1753294}. There, though $\Lm$ is contained in $\RE.$
For more general results in this direction see \cite{AM1,AM2}.

\end{remark}

\begin{remark}
In a series of papers
\cite{Vasyunin83},
\cite{Vasyunin84},
\cite{BNO},
\cite{Hartmann96}, \cite{Hartmann99}),
 the free interpolation problem has been studied
and
a description of traces of bounded analytic functions
on a finite union of Carleson sets
has been obtained in terms of divided differences.
In view of  well known connections between interpolation
and basis properties,  these results may be partially
(\cite{Hartmann96}, \cite{Hartmann99})
rewritten
in terms of geometrical properties
of exponential DD in \liy.
\end{remark}

\section{Main Results}\LB{n1}
Let
$\Lm=\{\lm_n\}$ be a sequence in $\CX$
ordered in such a way that $\Re \lm_n$ form
a nondecreasing sequence. We connect with $\Lm$ the exponential family
$$
\E(\Lm)=\{
e^{i\lmn t}, te^{\lmn t}, \dots, t^{m_{\lmn}-1}e^{i\lmn t} \},
$$
where $m_{\lmn}$ is the multiplicity  of $\lmn\in \Lm$.

For the sake of simplicity, we confine ourselves to
the case
$\sup |\Im \lm_n| < \iy$. The multiplication \op
$f(t)\mapsto e^{-at}f(t)$ is an isomorphism in \lT for any $T$
and  maps  exponential functions $e^{i\lm_n t}$ to
$e^{i(\lm_n + ia)t}$. Since we are interesting in the Riesz basis
property of linear combinations of functions $t^r e^{i\lm_n t}$
in \lT ,  we can suppose without loss of generality that
 $\Lm$ lies in a strip $S:=\{z|\,0<\a \le \Im \lm\le \be<\iy\}$
in the upper half plane.

The sequence  $\Lm$ is called
{\it uniformly discrete}
or {\it separated} if condition (\ref{discr}) is fulfilled.
Note that in this case all points $\lm $ are simple and we
 do not need
 differentiate between a sequence and a set.

We say that  $\Lm$
{\it relatively uniformly discrete} if $\Lm$ can be
decomposed into a finite number of  uniformly discrete subsequences.
Sometimes we shall simply say that such a $\Lm$ is a finite union of
uniformly discrete sets, however we always consider a point $\lambda_n$
to be assigned a
multiplicity.

\subsection{Splitting of the spectrum and subspaces
of exponentials}\LB{n12}
Here we introduce notations needed to formulate the main result.
For any $\lm\in\CX$, denote by $ D_\lm(r)$ a disk with center $\lm$
and radius $r$. Let $G^{(p)}(r)$, $p=1,2,\dots,$ be the
connected components of the union
$\cup_ {\lm\in \Lm} D_\lm(r)$. Write
$\Lm^{(p)}(r)$ for the subsequences of $\Lm$
lying in $G^{(p)}$,
$\Lm^{(p)}(r):=
\{\lmn |\, \lmn \in  G^{(p)}(r)\} $, and $\L^{(p)}(r)$ for subspaces spanned
by corresponding exponentials
$\{t^n e^{i\lm t} \}$, $\lm\in \Lm^{(p)}(r), n=0,\dots, m_{\lm}-1$.

\begin{lemma}\LB{card}
Let   $\Lm$ be a union of $N$ \uds s $\Lm_j$,
$$
\dl_j:=\dl(\Lm_j):=\inf_{\lm\neq \mu;\;\lm,\mu\in\Lm_j }|\lm- \mu|, \
\dl:=\min_j \dl_j.
$$
Then for
$$
r<r_0:=\frac{\dl}{2N}
$$
the number $\Np(r)$ of elements of $\Lm^{(p)}$ is at most $N$.
\end{lemma}

\proof
Let points $\mu_k$, $k=1,\dots,N+1$, belong to the same $\Lm^{(p)}$.
Then the distance between any two of these points is less than
$2rN$ and so less than $\dl$. From the other hand, there at least
two among $N+1$ points which belong to the same
$\Lm_j$ and, therefore, the distance between them is not less than
$\dl$. This contradiction  proves the lemma.
\qed

We call by {\it \Lb}
a family in a Hilbert space which forms a \Rb
in the closure of its linear span.

The following statement is a small modification of the theorem of
Vasyunin \cite{Vasyunin78}.
\begin{lemma}\LB{Vasyunin78}
Let $\Lm$ be a \rds. Then
for any $r>0$ the family
of subspaces $\L^{(p)}(r)$ forms an \Lb in \liy.
\end{lemma}
This statement is proved in
\cite{Vasyunin78} and
\cite[Lect. IX]{Nik} for $r=r_0/16$ and for disks
in so called hyperbolic metrics, however one can easily check that
the proof with obvious modifications remains valid for all $r$.

In applications we often meet the case of real $\Lm's$.
Then a relatively discrete set can be  characterized    using another
parameter than in Lemma \ref{card}. The following statement can be easily
proved similar to Lemma \ref{card}.

\begin{lemma}\LB{delta}
A real sequence $\Lm $ is a union of $N$ \uds s $\Lm_j $
if and only if
$\inf_n (\lm_{n+N}- \lm_n):= \tl\dl>0$ for all $n$. Along with that,
$\min_j \dl(\Lm_j)\le \tl\dl$
\end{lemma}

\subsection{Divided differences}\LB{n13}

Let  $\mu_k$, $k=1,\dots,m$, be arbitrary complex numbers,
not necessarily distinct.
\begin{definition}\LB{gdd_def}
Generalized divided difference (GDD)
of  order zero of $e^{i\mu t}$ is
$[\mu_1](t):=e^{i\mu_1 t}$.
GDD of the order $n-1$, $n\le m$, of
$e^{i\mu t}$ is
$$
[\mu_1,\dots,\mu_n]:=
\cases{\ds
\frac{[\mu_1,\dots,\mu_{n-1}]-[\mu_2,\dots,\mu_n] }{\mu_1-\mu_n},
                &$\mu_1\ne\mu_n$,\cr
\frac{\partial}{\partial \mu}[\mu,\mu_2,\dots,\mu_{n-1}]\Big |_{\mu= \mu_1},
                &$\mu_1=\mu_n$.
}
$$
\end{definition}

If all
 $\mu_k$ are distinct
one can easily derive the explicit formulae for GDD:
\EQ{dd}
{
 [\mu_1,\dots,\mu_n]=\sum_{k=1}^n \frac{e^{i\mu_k t}}{\prod_{j
\neq k} (\mu_k - \mu_j)}.
}

For any points $\{\mu_k\}$ we can write \cite[p. 228]{Shilov}
\begin{eqnarray} \LB{dd_int}
\nonumber
 [\mu_1,\dots,\mu_n]=
\int_0^1d\tau_1
\int_0^{\tau_1}d\tau_2 \dots
\int_0^{\tau_{n-2}}d\tau_{n-1} (it)^{n-1}\\
\exp\Big( it
\l[ \mu_1+\tau_1(\mu_2-\mu_1)
+\dots + \tau_{n-1}(\mu_n-\mu_{n-1})\r]\Big).
\end{eqnarray}

\begin{theorem}\LB{gdd_1}
The following statements are true.

(i)
Functions
$\f_1:=[\mu_1],\dots,\f_n:=[\mu_1,\dots,\mu_n]$,
depend on parameters $\mu_j$
continuously and  symmetrically.
If
points
$\mu_1,\dots,\mu_n$ are in a convex domain $\Omega\subset\CX,$ then
\EQ{iii}{
|\f_j(t)|\le c_n e^{\gm t}, \gm:=-\inf_{z\in\Omega}\Im z.
}

(ii)
Functions
$\f_1,\dots,\f_n$,
are linearly independent.

(iii)
If points
$\mu_1,\dots,\mu_n$
are distinct, then
the family of GDD's forms a basis in the span of exponentials
$e^{i\mu_1t},\dots,e^{i\mu_n t}$.

(iv) Translation of the set $\mu_1,\dots,\mu_n$ leads to
multiplying of GDD's  by an exponential:
$$
[\mu_1+\lm,\dots,\mu_n+\lm]=e^{i\lm t}[\mu_1,\dots,\mu_n].
$$

(v)
For any $\eps>0$ and any $N\in \NN$, there exists $\dl$
such that
the estimates
$$
||[\mu_1,\dots,\mu_j](t)-e^{i\mu t} t^{j-1}/(j-1)!||_{\liy}
< \eps,\ j=1,\dots,N,
$$
are valid
for any  $\mu$ in the strip $S$ and all points
$ \mu_1, \ldots , \mu_N$ belonging to the disk $D_\mu(\dl)$ of  radius
$\dl$ with the center at $\mu$.
\end{theorem}

\subsection{Bases of elements}\LB{n14}
Let  $\Lm^{(p)}(r)=\left\{\lm_{j,p}\right\}, j=1,\ldots,\Np(r)$
 be subsequences described in subsection \ref{n12}
Denote by $\{\E^{(p)}(\Lm,r) \}$
the family of GDD
corresponding to the points $\Lm^{(p)}(r)$:
$$
\E^{(p)}(\Lm,r)=\{
[\lm_{1,p}],[\lm_{1,p},\lm_{2,p}],\dots,
[\lm_{1,p}, \dots, \lm_{\Np,p}] \}.
$$
Note that $\E^{(p)}(\Lm,r)$ depends on enumeration
of $\Lm^{(p)}(r)$, although every GDD depends symmetrically on its parameters,
see the assertion (i) of the last theorem.

\begin{theorem}\LB{main_21}
Let $\Lm$ be a relatively uniformly discrete sequence and $r<r_0$. Then

(i) the family $\{\E^{(p)}(\Lm,r)\}$
forms a \Rb in $L^2(0,T)$ \iff\
 there exists an entire function $F$ of
exponential type with indicator diagram of width $T$
and zeros at points   $\lmn$  multiplicity $m_{\lmn}$
{\rm (the generating function of the family
$\E(\Lm)$ on the interval $(0,T)$)}
 such that, for some real $h$, the function
$\left|F(x+ih)\right|^2$ satisfies
the Helson--Szeg\"{o}\ condition (\ref{HS});

(ii) for any finite sequence $\{a_{p,j}\}$ the \ie
$$
\|\sum_{p,j}a_{p,j}e^{i\lm_{j,p} t}\|^2_{L^2(0,T)} \geq C
\sum_{p,j}|a_{p,j}|^2\dl_p^{2(\Np(r) -1)}
$$
is valid with  a constant $C$ independent of $\{a_{p,j}\}$, where
$$
\dl_p:=\min\{|\lm_{j,p}  - \lm_{k,p}  |\, \Big |\,k\ne j \}.
$$
\end{theorem}

Suppose now that $\Lm$ is not a relatively uniformly discrete sequence.
Then,
$$
\sup_p\Np(r)=\iy \ \ \mbox{ for any} \ \ r>0.
$$
It is possible also that there is an  infinite set $\Lm^{(p)}$.

We show that in this case the family of GDDs
is not uniformly minimal  even in $\liy$ at least, for  some enumeration of
points.

\begin{theorem}\LB{nonrel}
For any  $r>0$, there exists numbering of
points in the subsequences $\Lm^{(p)}$
such that family $\{\E^{(p)}(\Lm,r)\}$ is not uniformly minimal
in $L^2(0,\iy)$. Moreover,
for any $\eps$,  there exists
a pair $\f_m,\f_{m+1}$ of GDD,
corresponding to some set $\Lm^{(p)}(r)$ such that
$$
\on{angle}_{\liy}(\f_m,\f_{m+1})<\eps.
$$
\end{theorem}
\subsection{Application to an observation problem}\LB{appl}

Before  starting the proofs of main results,
we present an application of  Theorem \ref{main_21}
to observability  of a coupled 1d system studied in \cite{KomornikL00}:
$$
\begin{cases}
{
\frac{\d^2}{\d t^2}u_1-
\frac{\d^2}{\d x^2}u_1 +Au_1+Bu_2=0 &
                                                        in $(0,\pi)\times\RE$,\cr
\frac{\d^2}{\d t^2}u_2+\frac{\d^4}{\d x^4}u_2 +Cu_1+Du_2=0 &
                                                        in $(0,\pi)\times\RE$,\cr
u_1=u_2=\frac{\d ^2}{\d x^2}u_2=0   &
                                                                for $ x=0,\pi$,\cr
u_1=y_0\in H^1_0(0,\pi),\  \frac{\d }{\d t}u_1=y_1\in L^2(0,\pi),   &
                                                  for $ t=0,$     \cr
u_2=\frac{\d }{\d t}u_2=0                                       &
        for $ t=0$             \cr
}
\end{cases}
$$
($A$, $B$, $C$, $D$ are constants).

We introduce the initial energy $E_0$ of the system,
$E_0:=\|y_0\|^2_{H^1} +\|y_1\|^2_{L^2}$.

In the paper of
C.~Baiocchi, V.~Komornik, and P.~Loreti
\cite{KomornikL00}
the partial observability, i.e. inequality
\EQ{obs}{
\l\Vert\frac{\d}{\d x} u_1(0,t)\r\Vert_{L^2(0,T)}^2 \geq c E_0
}
with a constant $c$ independent of $y_0$ and $y_1$, has been proved
for almost all 4--tuples
$(A,B,C,D)$ and for $T>4\pi$.
(It means that we can recover the initial state
via the observation
$\l\Vert\frac{\d}{\d x} u_1(0,t)\r\Vert_{L^2(0,T)}^2$
during the time $T$ and the operator:
{\it
observation $\to$  initial state
}
is bounded.
The authors conjectured the system is probably partially observed for
$T>2\pi$.  Here  we  demonstrate this fact using the basis property of
exponential DD.
\begin{proposition}\LB{p_obs}
{}For almost all 4--tuples $(A,B,C,D)$ and for $T>2\pi$ the estimate
(\ref{obs}) is valid.
\end{proposition}

To prove this proposition
we use the representation and properties of the solution given in the paper
\cite{KomornikL00}.
To apply the Fourier method we introduce the eigenfrequencies
$\om_k$, $\nu_k$, $k\in \NN$ of the system, where
 $\nu_k^2$  and $\om_k^2$ are the eigenvalues of the matrix
$$
\l({
\begin{array}{cc}
k^2+A & B \cr
C     & k^4 +D
\end{array}
}
\r) .
$$
It is easy to see
that the following
asymptotic relations are valid:
\EQ{asymp}{
 \nu_k=k+A/2k +O(k^{-3}), \
\om_k=k^2+D/2k^2 +O(k^{-6}).
}
We suppose that
all
 $\om_k$ and $\nu_k$ are distinct
(this is true for almost all 4-tuples).
Then the first component of the solution of the system
can be written in the form
\EQ{wave}{
u_1(x,t)=\sum_{k\in \KKK}
\l[\a_ke^{i\om_kt} +\be_ke^{i\nu_kt}\r] \sin kx,
}
where
$
\KKK:=\ZZ \backslash \{0\}, \ \om_{-k}:=-\om_k,  \ \nu_{-k}:=-\nu_k,
$

The coefficients $\a_k$, $\be_k$ entered the last sum can be
expressed via the initial data, and the authors of
\cite{KomornikL00} show that under  zero initial condition for
$u_2$, we have
\EQ{KL}{
|\a_k|^2 + |\a_{-k}|^2\prec k^{-8}(|\be_k|^2 + |\be_{-k}|^2)
}
and
\EQ{energy}{
E_0\asymp
\sum_{k\in \KKK} k^2  |\be_k|^2.
}
Relations (\ref{KL}) and  (\ref{energy}) mean, correspondingly,
one--sided and two--sided inequalities with constants
which do not depend on sequences $\{\a_k\}$ and $\{\be_k\}$.

We now have to study the exponential family
$$
\E=\{e^{i\lm t}\}_{\lm\in \Lm},\
\Lm= M \cup \Om,\
  M=\{\nu_k\}, \ \Om=\{\om_k\}, \ k\in \KKK.
$$
For $r<1$ and for large enough $|\lm|$
the family $\Lm^{(p)}(r)$ consist of one point $\nu_k$ if k
is not a full square or of two points
$\nu_{\on{sign }k\,k^2}$ and $\om_k$. Denote by $\E^{(p)}(\Lm,r)$
the related family of exponential DD.

Suppose for a moment that there exists
an entire function $F$ of
exponential type with indicator diagram of width $T>2\pi$
vanishing at  $\Lm$ ($F$ may also have another zeros) such that
$|F(x+ih)|^2$ satisfies
the Helson--Szeg\"o condition (\ref{HS}) for some real $h$.
Then, by  Theorem \ref{main_21}, the family $\E^{(p)}(\Lm,r)$
  forms an \Lb in \lT.
Expanding exponentials in  DD of the zero and first order,
$$
\a e^{i\lm t}+\be e^{i\mu t}=(\a+\be)e^{i\lm t}-\be(\lm-\mu)[\lm,\mu],
$$
for finite sequences
$\{p_k\}$, $\{q_k\}$
we obtain
the estimate
\EQA{2pi}{rl}{
\ds
\l\Vert\sum_{k\in \KKK}
\l[ p_ke^{i\om_kt} +q_ke^{i\nu_kt}\r]
\r\Vert^2_{L^2(0,T)}  \asymp
\ds
\sum_{k\in \KKK, \sqrt{|k|}\notin {\bf N} } |q_k|^2+\\
\sum_{k\in \KKK}\l[|p_k + q_{\on{sign} k\,k^2}|^2 +
|\om_k-\nu_{\on{sign} k\,k^2}|^2
|q_{\on{sign} k\,k^2}|^2\r].
}
In the right hand side of this relation
the first sum corresponds to one dimensional
subspaces $\L^{(p)}(r)$
(this sum is
taken for all $k$ which are not full squares).
The second sum corresponds to two
dimensional subspaces.
Together with (\ref{wave}), it follows
\begin{eqnarray*}
\l\Vert \frac{\d}{\d x}u_1(0,t)\r\Vert^2_{\lT}=\l\Vert
\sum_{k\in \KKK}k
\l[\a_ke^{i\om_kt} +\be_ke^{i\nu_kt}\r]\r\Vert_{\lT}^2
\\
\succ
\sum_{k\in \KKK, \sqrt{|k|}\notin {\bf N} } k^2  |\be_k|^2
+
\sum_{k}\l[
|\om_k-\nu_{\on{sign} k\,k^2}|^2|k^2\be_{\on{sign} k\,k^2}|^2+
 \l|k\a_k+k^2\be_{\on{sign} k\,k^2}\r|^2
\r]\,.
\end{eqnarray*}

For $R$ large enough using (\ref{KL}) we have
$$
\sum_{|k|>R} \l[k\a_k+k^2\be_{\on{sign} k\,k^2}\r]
\succ
\sum_{|k|>R} k^2  |\be_{\on{sign} k\,k^2}|^2.
$$

Thus, we obtain
$$
\l\Vert \frac{\d}{\d x}u_1(0,t)\r\Vert^2_{\lT}\succ
\sum_{k\in \KKK} k^2  |\be_k|^2 \asymp E_0.
$$
In order to complete the proof of Proposition \ref{p_obs},
it  remains to construct
the function $F$. We do that using the following proposition which is
interesting in its own right.

\begin{proposition}\LB{mean}
Let $\Lm=\{\lm_n\}$ be a zero set of a sine type function
(see definition in
\cite[p. 61]{AI95}),
with the indicator diagram of with $2\pi$,
$\{\dl_n\}$ a bounded sequence of complex numbers,
and $F$ an entire function of the Cartwright class (\cite[p. 60]{AI95})
with the zero set  $\{\lm_n+\dl_n\}$.
If
$$
\lim_{N\to \iy}\sup_n\frac1N
\Big|\Re\l(\dl_{n+1} + \dl_{n+2}+\dots + \dl_{n+N}\r)\Big|=:d<\frac14,
$$
then $F$  has the indicator diagram of with $2\pi$, and
for any $d_1>d$, $h$,  functions
$
u, v \in L^{\infty}(\RE),
{\Vert v \Vert}_{L^{\infty}(\srline)} < 2\pi d_1
$
 may be found such that
$$
|F(x+ih)|^2=\exp\{u(x)+\tl v(x)\}
$$
for any real $h$ such that $|h|>\sup |\Im (\lmn + \dl_n)|$.
\end{proposition}

For the proof of this assertion, one argues using  \cite[Lemma 1]{Avd74}, and in a
similar way to the proof of \cite[Lemma 2]{Avd74}  that
there exists a sine type
function with zeros $\mu_n$ such
that, for any $d_1>d,$
$$
d_1\Re(\mu_{n-1}-\mu_n)\leq \Re (\lmn+\dl_n - \mu_n)\leq d_1
\Re(\mu_{n+1}-\mu_n)\,.
$$
Then Proposition  \ref{mean} follows directly from \cite[Lemmas 1, 2]{AJ}.

We begin to construct $F$ putting
$$
F_1(z):=z  \prod_{n\in{\bf N} }(1-\frac{z^2}{\nu_n^2}).
$$
In view of the asymptotics (\ref{asymp}) of $\nu_k$ and Proposition \ref{mean},
$F_1(z)$ may be written in the form
\EQ{HS1}
{
|F_1(x+ih)|^2=\exp\{u(x)+\tl v(x)\}
}
with
$
u, v \in L^{\infty}(\RE),
{\Vert v \Vert}_{L^{\infty}(\srline)} <  \eps
$
for any positive $\eps$.

Further, for any fixed $\dl >0$,
we  take the function
$\sin \pi \dl z/2$ and  for every $k \in \KKK $
find the zero $2n_k/ \dl $ of this function
which is
the nearest to $\om_k$.
 We  set
$$
F_2(z):=z  \prod_{n\in{\bf N} }(1-\frac{z^2}{\mu_n^2})\,,
$$
where
$$
\mu_n=\cases{\ds 2n/ \dl&\mbox{for } $n\notin \{n_k\}$ ,\cr
\om_k   & \mbox{for } $n= n_k$.
}
$$

Using  again
Proposition \ref{mean} we conclude that
for any $\eps>0$ function $|F_2(x+ih)|^2$  has the representation
similar to (\ref{HS1})
and the  width of the indicator diagram of $F_2$ is $2\dl$.

Taking $F=F_1F_2$ we obtain the  function with indicator diagram of
 width $2(\pi +\dl)$ vanishing on $\Lm$ and satisfying (\ref{HS}).
Proposition \ref{p_obs} is proved.

\section{Proofs}\LB{n2}

\subsection{
Proof of Theorem \ref{gdd_1}
}\LB{n21}

(i) Continuity on parameters $\{ \mu_k\}$ and estimate (\ref{iii})
follows immediately from the representation (\ref{dd_int}).
{}From (\ref{dd}) symmetry
is clear if all  $\{ \mu_k\}$  are different. Then, by continuity,
we obtain symmetry for any points in the sense that if
$\sigma$ is a permutation of $\{1,2,\dots,n\}$, then
$[\mu_1,\dots, \mu_n]=
[\mu_{\sigma(1)},\dots, \mu_{\sigma(n)}]$.

Now we describe the structure of GDD.
Let $z_1,z_2,\dots,z_n$ be complex numbers, not necessarily different.
Let us represent this set
as the union of $q$
different points,
$\nu_1, \dots,\nu_q$, with multiplicities
$m_1, \dots,m_q$;
($m_1+m_2+\dots+m_q=n$).

\begin{lemma}\LB{structure}
The GDD
$[z_1, \dots,z_n]$
of order $n-1$
is a linear combination of functions
$$
t^me^{i\nu_kt},\ k=1,\dots,q,\ m=0,1,\dots,m_k-1
$$
and the coefficients of  the leading  terms
$t^{m_k-1}e^{i\nu_kt}$ are not equal to zero.
\end{lemma}

{\sc Proof}\  of the lemma.
The statement is  clear for $n=1$.
Let us suppose that it is true for
GDD of the order $n-2$.

If $z_1\ne z_n$ then, by definition,
$[z_1, \dots,z_n]=\ds \frac{[z_1,\dots,z_{n-1}]-[z_2,\dots,z_n] }{z_1-z_n}$
and is a linear combination of the GDD's of order $n-2$.
Multiplicity of the point
$z_n$ in the set $\{z_2,\dots,z_n\}$ is more than
in the set $\{z_1,\dots,z_{n-1}\}$. Let $k$ be  such a
number that $z_n=\nu_k$. Then
the leading term
$t^{m_k-1}e^{i\nu_kt}$ appears  in $[z_2,\dots,z_n]$
by the inductive conjecture
and does not appear in $[z_1,\dots, z_{n-1}]$.

Let $z_1= z_n$. Then $[z_1,\dots,z_{n-1}]$ contains the leading term
$t^{m_1-2}e^{i\nu_1t}$. After differentiation in $z_1$
we  get the leading term
$t^{m_1-1}e^{i\nu_1 t}$ with nonzero coefficient.
Using symmetry of GDD relative to points $z_1,...,z_n$,
we complete the proof of
the lemma.
\qed

We are able now, with  knowledge of
the structure of GDD, to continue the proof of Theorem~\ref{gdd_1}.

(ii)
This assertion is the consequence of the `triangle' structure of GDD:
if we add $n$-th point, then a
GDD of the order $n-1$,
in comparison with  a GDD of the order $n-2$, contains either a new exponential or
a term of the form $t^me^{i\nu_pt}$ with the same frequency
$\nu_p$ and larger
exponent $m$.

(iii)
If all points
$\mu_1,\dots,\mu_n$
are different, then
the family of GDD's contains
all exponentials
$e^{i\mu_1t},\dots,e^{i\mu_n t}$.

(iv) Immediately follows from the definition.

(v) As well known, divided differences approximate
the derivatives of the corresponding order.
We use the following estimate \cite{Ullrich}.
\begin{proposition}
For any $\eps>0$,  $N\in \NN$, there exists $\dl$ such that
for any set $\{ z_j \}_{j=1}^N$ belonging to the disk $D_0(\dl)$ of a radius
$\dl$ with center at the origin, the estimates
\EQ{Ulr}
{
\Big|[z_1,\dots,z_j](t)- t^{j-1}/(j-1)!\Big|
< \eps,\ j=1,\dots,N,
}
are valid for $t \in [-\pi N, \pi N]$.
\end{proposition}
To proceed with (v) we choose $T$ large enough  that
\begin{eqnarray}\LB{v_1}
{
\Big\|[\mu_1,\dots,\mu_j]\Big\|_{L^2(T,\iy)}< \eps/3,\
\Big\|e^{i\mu t} t^{j-1}/(j-1)!\Big\|_{L^2(T,\iy)}< \eps/3,\
}
\nonumber
\\
 j=1,\dots,N,
\end{eqnarray}
for all $\mu,\mu_1,\ldots ,\mu_j$ lying in  the strip $S$.

Set
$\eps_1=\eps /3\sqrt{T}$,
 and let  $\dl$ be small enough in order to
(\ref{Ulr}) be fulfilled for such $\eps_1$,
$ t \in [0, T]$ and
any set $\{ z_j \}_{j=1}^{N}\in D_0(\dl)$.
Set $z_j:= \mu_j-\mu$. In view of (iv), we obtain
$$
\Big|e^{-i\mu t}[\mu_1,\dots,\mu_j](t)- t^{j-1}/(j-1)!\Big|
< \eps_1,\ j=1,\dots,N,\ t\le T,
$$
that implies
$$
\Big|[\mu_1,\dots,\mu_j](t)- e^{i\mu t}t^{j-1}/(j-1)!\Big|
< \eps/3\sqrt{T},\ j=1,\dots,N,
$$
for $t\le T$. Therefore, we have
$$
\Big\|[\mu_1,\dots,\mu_j](t)- e^{i\mu t}t^{j-1}/(j-1)!\Big\|_{\lT}
< \eps/3,\ j=1,\dots,N.
$$
Taking into account  (\ref{v_1}) we obtain the statement (v).
\qed

\subsection{
Proof of the Theorem \ref{main_21}.
 }\LB{n22}
\vskip2mm(i) From Lemma \ref{Vasyunin78} it follows that
the family $\{\L^{(p)}\}:=\{\L^{(p)}(r)\}$ forms
a \Rb in the closure of its span (i.e. an \Lb) in $L^2(0,\iy)$.
Introduce the projector $P_T$ from this closure {\it into} $L^2(0,T)$.
Clearly,
$\{\L^{(p)}\}$ forms a Riesz basis in $L^2(0,T)$
\iff $P_T$ is an isomorphism
{\it onto} $\lT$.
 This takes a place \iff the generating function of the family $\E(\Lm)$ satisfies
the Helson-Szeg\"{o} condition (\ref{HS})
(see \cite{KNP};
\cite{AI95}, Theorems II.3.14, II.3.17).

The \Lb property of the subspaces
$\L^{(p)}$ means that for any finite number of functions
$\psi_p\in \L^{(p)}$ we have the estimates:
\EQ{6.1}{
\|\sum_p a_p \psi_p\|_{L^2(0,\iy)}^2\asymp
\sum_p|a_p|^2 ||\psi_p ||_{L^2(0,\iy)}^2.
}
This relation means that there exist two positive constants $c$ and $C$
which do not depend on a sequence $\{a_p\}$ such that
$$
c\|\sum_p a_p \psi_p\|_{L^2(0,\iy)}^2\leq
\sum_p|a_p|^2 ||\psi_p ||_{L^2(0,\iy)}^2
\leq C\|\sum_p a_p \psi_p\|_{L^2(0,\iy)}^2.
$$

In each subspace $\L^{(p)}$ we choose the family of GDD's
corresponding to
$\lm_{j,p}\in \Lm^{(p)}$:
$\f_j^{(p)}:=[\lm_{1,p},\lm_{2,p},\dots,\lm_{j,p}]$.
$j=1,\dots, \Np$.
In view of Theorem \ref{gdd_1}, this family forms a basis in $\L\p$.

Let us
expand
$\psi_p$
in the basis
$\{\f_j\p\}_{j=1}^{j=\Np}$. Then statement (i) of the theorem
is equivalent to the estimates
\EQ{6.2}{
\|\sum_{p,j} a_{p,j}\f_j^{(p)}\|_{L^2(0,\iy)}^2\asymp
\sum_{p,j}|a_{p,j}|^2.
}
Taking into account (\ref{6.1}),
we see that
(\ref{6.2})
is equivalent to
\EQ{6.3}{
\|\sum_j a_j\f_j^{(p)}\|_{L^2(0,\iy)}^2\asymp
\sum_j|a_j|^2,
\mbox {uniformly in }p.
}

We introduce the $\Np\times\Np$ Gram matrices $\Gm\p$ corresponding to
the families
$\E^{(p)}$,
$$
\Gm\p:=\l\{(\f_k^{(p)},\f_j^{(p)})_{L^2(0,\iy)}\r\}_{k,j}.
$$
For an $\Np$-dimensional vector $a$ we have
$$
\|\sum_j a_j\f_j^{(p)}\|_{L^2(0,\iy)}^2=
\langle \Gm\p a,a\rangle,
$$
where
$\langle \cdot,\cdot\rangle$
is the scalar product in $\RE^{\Np}$.

In terms of the Gram matrices,
(\ref{6.3}) is true \iff\  matrices
$\Gm\p$ and their inverses are bounded uniformly in $p$.

\begin{lemma}\LB{gram}
Let different complex points
$\mu_1,\mu_2,\dots,\mu_n$
lie in a disk $D_\mu(R)\subset \CX_+$ of radius $R$
with the center $\mu$, $\mu\in S$;
$\f_1,\f_2,\dots,\f_n$ are corresponding GDD's,
and $\Gm$ is the Gram matrix of this family in $L^2(0,\iy)$.
Then the norms
$\langle\langle\Gm\rangle\rangle$ and
$\langle\langle\Gm^{(-1)}\rangle\rangle$ are estimated from
above by constants depending only on $R$ and $n$.
\end{lemma}
\proof
By Theorem~\ref{gdd_1},
functions
$\f_1(t),\f_2(t),\dots,\f_n(t)$
are estimated
above by constants depending only on $R$, $n$.
Therefore the entries of $\Gm$ are estimated by
${const}\,(R,n)$
and, then, the estimate of $\Gm$ from above is obtained.

We shall prove the estimate for the inverse matrices by contradiction.
Let us fix a disk and assume that for arbitrary $\eps>0$ there exist
different points
$\mu_1^{(\eps)},\mu_2^{(\eps)},\dots,\mu_n^{(\eps)}$
lying in the disk, and
normalized
$n$-dimensional vectors $a^{(\eps)}$  such that
for the corresponding Gram matrix $\Gm^{(\eps)}$ we have
\EQ{le_eps}{
\langle\Gm^{(\eps)}a^{(\eps)},a^{(\eps)}\rangle\le\eps.
}
Using compactness arguments we can choose a sequence
$\eps_n\to 0$ such that
$$
a^{(\eps_n)}\to a^0, \ \mu_j^{(\eps_n)}\to \mu_j^0.
$$
as $n\to \iy$.
Then the Gram matrix tends to the Gram matrix $\Gm^0$ for the
limit family and from (\ref{le_eps}) we see that
$$
\langle\Gm^0a^{0},a^0\rangle=0.
$$
Then
$
\sum a^0_j\f^0_j=0
$
and the limit family of GDD's
is linearly dependent that contradicts to Theorem~\ref{gdd_1}.

Thus, we have proved that the norm of $\Gm^{(-1)}$ is estimated from above
by constants depending only on $R$, $n$, and
$\mu$. In view of Theorem~\ref{gdd_1}(iv) the constant depends actually
not on $\mu$, but on $\Im \mu$.
Indeed,  translation of the disk on a real number,
$\mu \mapsto \mu+ x$, does not change the  Gram matrix.
Since $\a\leq \Im \mu \leq \beta,$ the norm of $\Gm^{(-1)}$ is
estimated uniformly in $\mu\in S$.
\qed

 {}From this lemma it follows that all Gram matrices
$\Gm_p$, $\Gm_p^{(-1)}$
are bounded uniformly in $p$. (It was supposed in the lemma that all points
$\mu_1,\mu_2,\dots,\mu_n$ are distinct.
For the case of multiple points
we use the continuity of DD on parameters.)
The assertion (i) of the Theorem \ref{main_21} is proved.

(ii)
 We need to  estimate the Gram matrices $\Gm$
for the exponential family
$\{e^{i\mu_1t},\dots,e^{i\mu_n t}\}$
in $L^2(0,T)$ from below,
where  $\{\mu_j\}$ is a fixed set $\Lm\p$ and $n=\Np$.
Since the projector
$P_T$  is an isomorphism (see (i)), we may do it for exponentials
on the positive semiaxis, i.e., in $L^2(0,\iy)$.
Denote by
$e_j(t)$ the normalized exponentials
$$
e_j(t):=\frac{1}{\sqrt{2\Im \mu_j}}e^{i \mu_j t}.
$$
For the Gram matrices $\Gm_0$, corresponding to
the normalized exponentials, we have
$$
\Gm=\on{diag}[\sqrt{2\Im \mu_j}]
\,\Gm_0\,
\on{diag}[\sqrt{2\Im \mu_j}],
$$
and we can estimate $\Gm_0$ instead $\Gm$, since
$|\Im\mu_j|\asymp 1$.
It can be easily shown that the inverse matrix is the Gram matrix
for the biorthogonal family
$e'_j(t)$ and
$$
(\Gm_0)^{(-1)}_{jj}=\|e'_j\|^2 =
\underset{k\ne j}{\Pi}{\Big |\frac{\mu_k-\bar\mu_j}{\mu_k-\mu_j}\Big |}^2
$$
(see \cite{Nik, AI95}).
Elementary calculations give
$$
(\Gm_0)^{(-1)}_{jj}=\|e'_j(t) \|^2 \prec \dl_p^{-2(n-1)}
$$
Then
$$
|(\Gm_0)^{(-1)}_{jk}|=(e_j',e_k')
\le \|e'_j \|\|e'_k \|
\prec \dl_p^{-2(n-1)}
$$
and so,
$$
\langle \langle \Gm_0^{(-1)} \rangle\rangle \prec \dl_p^{-2(n-1)}.
$$
The theorem is proved. \qed

\vskip2mm
\subsection{
Proof of Theorem \ref{nonrel}.
}
\vskip2mm

As well known, the family
$\{t^n\}$, $n=0,1,\dots,$
of powers is not minimal in $\lT$ for any $T$.
Moreover, direct calculations for
$\f^0_j:=e^{i\mu t} t^{j-1}/(j-1)!$
give
$$
\on{angle}_{\liy}(\f^0_m,\f^0_{m+1})\to 0,\ \ m\to\iy,
$$
uniformly in $\mu\in S$.
Since $\Lm$ is not a \rds,
for any $\dl>0$, $m\in \NN$ we are
able to find a disk $D_\mu(\dl)$ in the strip
$S$, which contains $m$ points of  $\Lm$.
In view of Theorem \ref{gdd_1}(v), GDD corresponding to these points
are $\eps$--close to `unperturbed' functions
$\f^0_j$.  This proves the assertion of Theorem \ref{nonrel}.

\end{document}